\documentclass[12pt]{amsart}

\usepackage{amsfonts,amsmath,amsthm}
\usepackage{latexsym}

\newtheorem{theorem}{Theorem}[section]

\newtheorem{remark}[theorem]{Remark}

\theoremstyle{definition}

\def \e {\varepsilon}

\def \d {\delta}

\def \< {\langle}
\def \> {\rangle}

\begin{document}

\title[Sets with no solutions to $x+y=3z$]{Sets with no solutions to $x+y=3z$}

\author[M. Matolcsi, I. Ruzsa]{M\'at\'e Matolcsi, Imre Z. Ruzsa}
\address{M. M., I. R. : Alfr\'ed R\'enyi Institute of Mathematics,
Hungarian Academy of Sciences POB 127 H-1364 Budapest, Hungary
Tel: (+361) 483-8302, Fax: (+361) 483-8333}
\email{matomate@renyi.hu, ruzsa@renyi.hu}

\thanks{The authors were supported by the ERC-AdG 228005, and OTKA Grants No. K81658.}


\begin{abstract}
This short note gives an upper bound on the measure of sets $A\subset [0,1]$ such that $x+y=3z$ has no solutions in $A$.
\end{abstract}

\maketitle

\bigskip



\section{Introduction}

In this note we consider measurable sets $A\subset [0,1]$ such that the linear equation $x+y=3z$ has no solutions in $A$. In particular, we prove that the measure of $A$ satisfies $|A|\le \frac{1}{2}-\delta$ with $\d =\frac{1}{114}$. To put this result in context let us briefly describe the history of the problem.

\medskip

A set $A$ of real numbers is called $k$-sum-free if it does not contain elements $a,b,c$ such that $a+b=kc$, where $k$ is a positive integer.
(It is customary to require that not all $a,b,c$ be equal, to avoid trivial solutions if $k=2$). Let $f(n,k)$ denote the maximal cardinality of a $k$-sum-free set in $\{1, 2, \dots n\}$. The quantity $f(n,k)$ and the possible structure of maximal $k$-sum-free sets $A(n,k)$ has been studied extensively over the past decades, and an almost complete understanding has been reached for all values of $k$, except $k=2$.

\medskip

For $k=1$ we have $f(n,1)=\lceil \frac{n}{2} \rceil$ and the extremal sets are well-known. For $n$ odd, there are two  maximal sum-free sets: the set of odd numbers in $\{1, \dots n\}$, and the "top half" $\{ \frac{n+1}{2} , \dots n\}$. For $n$ even, $n\ge 10$, there are three maximal sum-free sets: the set of odd numbers, $\{ \frac{n}{2}+1 , \dots n\}$, and $\{ \frac{n}{2} , \dots n-1\}$.

\medskip

For $k=2$ the famous theorem of Roth \cite{roth} says that $f(n,2)=o(n)$, and giving tight upper and lower bounds on $f(n,2)$ is a notoriously difficult and famous problem.

\medskip

For $k=3$ Chung and Goldwasser \cite{cg2} proved that $f(n,3)=\lceil \frac{n}{2} \rceil$ for any $n\ne 4$, and the set of odd numbers is the unique maximal 3-sum-free set for $n\ge 23$.

\medskip

For $k\ge 4$ Chung and Goldwasser \cite{chung} gave an example of a $k$-sum-free set of cardinality $\frac{k(k-2)}{k^2-2}n+\frac{8(k-2)}{k(k^2-2)(k^4-2k^2-4)}n+O(1)$, and conjectured that $\lim_{n\to \infty } \frac{f(n,k)}{n}= \frac{k(k-2)}{k^2-2}+\frac{8(k-2)}{k(k^2-2)(k^4-2k^2-4)}$. They also conjectured that if $n$ is large then the extremal $k$-sum-free sets consist of three intervals of consecutive integers with slight modifications at the end-points. Subsequently Baltz \& al. \cite{baltz} proved the first conjecture and gave a structural result on maximal $k$-sum-free sets which is very close to proving the second.

\medskip

The examples for $k\ge 4$ in \cite{chung} were based on the solution of the continuous version of the problem. 
Namely, let $\mu (k)$ denote the maximal possible measure of a measurable $k$-sum-free set $A\subset [0,1)$. Intuitively one expects that information on $\mu(k)$ and the structure of the maximal sets should provide information on $\lim_{n\to \infty } \frac{f(n,k)}{n}$. This is indeed the case for most values of $k$. 

\medskip

For $k=1$ it is easy to see that $\mu(1)=\frac{1}{2}$ and the "top-half" interval $(\frac{1}{2}, 1]$ is essentially the only maximal 1-sum-free set. This is in analogy with $f(n,1)$, with the exception that the extremal set of odd numbers does not have a continuous analogue. 

\medskip

For $k=2$ a simple Lebesgue-point argument shows that $\mu(2)=0$, in analogy with $f(n,2)=o(n)$.

\medskip

For $k\ge 4$ the analogy continues to hold, as Chung and Goldwasser \cite{chung} showed that $\mu(k)=\frac{k(k-2)}{k^2-2}+\frac{8(k-2)}{k(k^2-2)(k^4-2k^2-4)}$, with the extremal set being a union of three intervals. 

\medskip

Therefore, the only case left open is $k=3$ in the continuous setting, the value of $\mu(3)$ being unknown. The largest known 3-sum-free set in the interval $(0,1]$, also given in \cite{chung}, is
\begin{equation}\label{conj}
A= \left(\frac{8}{177},\frac{4}{59}\right)\cup\left(\frac{28}{177},\frac{14}{59}\right)\cup\left(\frac{2}{3},1\right),
\end{equation}
with $|A|=\frac{77}{177}$. (Of course, some endpoints of the intervals can be included in $A$ but it does not change the measure.) 
In fact, Chung and Goldwasser conjecture that if $|A|$ is maximal and the set $A$ itself is maximal with respect to inclusion then $A$ must be equal to the union of these three intervals together with some of their endpoints \cite[Conjecture 3]{chung}. Let us remark here that we have made some  computer experiments which support this conjecture: an easy linear programming code shows that if $A$ is the union of 2, 3, 4 or 5 disjoint intervals then $|A|\le \frac{77}{177}$, with equality holding only for the intervals above.

\medskip

The primary motivation of this note is to show that the analogy between the discrete and continuous versions of the problem breaks down for $k=3$. Namely, $\lim_{n\to \infty } \frac{f(n,3)}{n}=\frac{1}{2}$ as mentioned above, but $\mu(3)\le \frac{1}{2}-\frac{1}{114}$ according to Theorem \ref{thm1} below. This shows that the extremal set of odd numbers for $f(n,3)$ cannot have a continuous analogue in $(0,1]$. In fact, we are fairly convinced that \cite[Conjecture 3]{chung} is true, and $\mu(3)=\frac{77}{177}$ but we are currently unable to prove it.

\medskip

\section{Sets with no solutions to $x+y=3z$}\label{sec2}

In this section we prove our main result.

\begin{theorem}\label{thm1}
Let $A\subset [0,1]$ be a measurable set such that there exist no $x,y,z\in A$ for which $x+y=3z$ holds. Then the measure of $A$ satisfies $|A|\le \frac{1}{2}-\delta$ with $\d =\frac{1}{114}$.
\end{theorem}

\begin{proof}

First note that we can assume that $A$ is closed, because the Lebesgue measure is inner regular. Second, we can assume that $1\in A$ because otherwise we could consider an appropriate dilate $\alpha A$ of $A$ with some $\alpha > 1$. With these assumptions $\inf (A)$ belongs to $A$, and
$\rm{diam} (A)=1-\inf (A)$. Also, let us introduce the notations $x=|A|$, $a=\inf (A)$ and $C=\frac{A+A}{3}$. Note that $C\subset [\frac{2a}{3}, \frac{2}{3}]$, $A\cup C\subset [\frac{2a}{3}, 1]$, and by assumption $A$ and $C$ are disjoint. We will first prove $x\le \frac{1}{2}$.

\medskip

Applying Corollary 3.1 in \cite{ruzsa1} with $A=B$ we obtain that $|A+A|\ge \min (3|A|, |A|+\rm{diam}(A))=\min (3x, x+1-a)$. This implies $|C|\ge \min (x, \frac{x+1-a}{3})$. If the minimum here is $x$ then using the fact that $A$ and $C$ are disjoint in $[\frac{2a}{3}, 1]$ we obtain $2x\le 1- \frac{2a}{3}$, and hence
\begin{equation}\label{e1}
x\le \frac{1}{2}-\frac{a}{3}
\end{equation}
follows. If the minimum is $\frac{x+1-a}{3}$ then we obtain $x+\frac{x+1-a}{3}\le 1-\frac{2a}{3}$, which implies
\begin{equation}\label{e2}
x\le \frac{1}{2}-\frac{a}{4}.
\end{equation}
In both cases we conclude that $x\le \frac{1}{2}$. Moreover, by dilating the interval $[0,1]$ we obtain that any set $\tilde{A}\subset [0,w]$ which contains no solutions to $x+y=3z$ satisfies
\begin{equation}\label{e3}
|\tilde{A}|\le \frac{w}{2}.
\end{equation}

\medskip

Assume now, by contradiction, that $x=1/2-\d$ with $\d < \frac{1}{114}$. From the argument above we see that in this case
\begin{equation}\label{adelta}
a\le 4\d
\end{equation}
must hold. Also, $C$ does not intersect $(\frac{2}{3}, 1]$, so $A$ must contain most of this interval, with the possible exception of a set of measure at most $2\d$ for the following reason. If $|A\cap (\frac{2}{3},1]|= \frac{1}{3}-\e$ then we get a term $-\frac{\e}{2}$ and $-\frac{3\e}{4}$ on the right hand side of \eqref{e1} and \eqref{e2}, respectively, and if $\e >2\d$ then this would imply $|A|\le \frac{1}{2}-\d$ in both cases.

\medskip

Taking into account that $a\in A$ we get that $|(A+A)\cap (a+\frac{2}{3}, a+1]|\ge \frac{1}{3}-2\d$, and  $|(A+A)\cap (\frac{4}{3}, 2]|\ge \frac{2}{3}-4\d$. Dividing by 3 and using again the disjointness of $A$ and $C$ we get that $|A\cap (\frac{a}{3}+\frac{2}{9}, \frac{a}{3}+\frac{1}{3}]|\le \frac{2\d}{3}$ and $|A\cap (\frac{4}{9}, \frac{2}{3}]|\le \frac{4\d}{3}$.
Let us simply disregard these parts of $A$ and let $A'=A\setminus  ( (A\cap (\frac{a}{3}+\frac{2}{9}, \frac{a}{3}+\frac{1}{3})) \cup (A\cap (\frac{4}{9}, \frac{2}{3})))$. Then $|A'|\ge |A|-2\d$.

\medskip

Note that $A'$ naturally breaks up into three disjoint parts: $A_1=A'\cap [a, \frac{a}{3}+\frac{2}{9}]$, $A_2=A'\cap [\frac{a}{3}+\frac{1}{3}, \frac{4}{9}]$, and $A_3=A' \cap [\frac{2}{3},1]$. Note that by \eqref{e3} and \eqref{adelta} we have
\begin{equation}\label{e4}
|A_1|\le \frac{a}{6}+\frac{1}{9}\le \frac{2\d}{3}+\frac{1}{9}.
\end{equation}
Let $b=\sup (A_1)$, $c=\inf (A_2)$ and $d=\sup (A_2)$. We have $b\le \frac{a}{3}+\frac{2}{9}$, and $\frac{a}{3}+\frac{1}{3}\le c\le d\le \frac{4}{9}$. If $A_2$ is non-empty then $a,b,c,d\in A'$. If $A_2$ happens to be empty, then $|A|\le |A'|+2\d =|A_1|+|A_3|+2\d \le (\frac{2\d}{3}+\frac{1}{9})+(\frac{1}{3})+2\d \le \frac{4}{9}+\frac{2\d}{3}+2\d < \frac{1}{2}-\d$, which is a contradiction.

\medskip

Thus $A_2$ is non-empty, and $A_1+A_3 \supset (a+A_3)\cup (b+A_3)$ and this latter set equals the interval $[a+\frac{2}{3}, b+1]$ with the exception of a set of measure at most $4\d$. Similarly, $A_2+A_3 \supset (c+A_3)\cup (d+A_3)$ which equals the interval $[c+\frac{2}{3}, d+1]$ with the exception of a set of measure at most $4\d$. Let us now consider two cases, according to the magnitude of $c-b$.

\medskip

If $c-b>\frac{1}{3}$, then $|A|\le |A'| +2\d =|A_1|+|A_2|+|A_3|+2\d \le (b-a)+(d-c)+\frac{1}{3}+2\d = (d-a)+(b-c)+\frac{1}{3}+2\d \le \frac{4}{9}-a-\frac{1}{3}+\frac{1}{3}+2\d \le \frac{4}{9} +2\d < \frac{1}{2}-\d$, a contradiction.

\medskip

If $c-b\le \frac{1}{3}$, then the intervals $[a+\frac{2}{3}, b+1]$ and $[c+\frac{2}{3}, d+1]$ overlap, so that $(A_1+A_3)\cup (A_2+A_3)$ equals the interval $[a+\frac{2}{3}, d+1]$ with the exception of a set of measure at most $8\d$. Therefore set $C$ contains the interval $I=[\frac{a}{3}+\frac{2}{9}, \frac{d}{3}+\frac{1}{3}]$, with the exception of a set of measure at most $\frac{8\d}{3}$. Notice, however, that $\frac{d}{3}+\frac{1}{3}>d$ because $d<\frac{1}{2}$. Therefore, the interval $I$ fully covers $A_2$, and hence $|A_2|\le \frac{8\d}{3}$.  Therefore, $|A|\le |A'|+2\d =|A_1|+|A_2|+|A_3|+2\d \le (\frac{a}{6}+\frac{1}{9})+(\frac{8\d}{3})+(\frac{1}{3})+2\d =\frac{4}{9}+\frac{16}{3}\d < \frac{1}{2}-\d$, a contradiction.

\end{proof}

\begin{remark}\rm
It is more than likely that one could obtain a larger $\delta$ by doing the case-by-case analysis with a little more care. However, the arguments above do not seem to lead to the structural result that the extremal 3-sum-free set must be the union of the three intervals stated in \cite[Conjecture 3]{chung}.
\end{remark}

\end{document}